\documentstyle[12pt]{article}
\title{From Flows and Metrics to Dynamics}
\author{C.Udri\c ste and A.Udri\c ste \\ {} \\
University Politehnica of Bucharest \\
Department of Mathematics I \\
Splaiul Independen\c tei 313 \\
77206 Bucharest, Romania \\
email:udriste@mathem.pub.ro}

\date{}
\begin{document}
\maketitle

\newcommand{\ty}{\infty}
\newcommand{\ov}{\over}
\newcommand{\di}{\displaystyle}
\newcommand{\si}{\sigma}
\newcommand{\na}{\nabla}
\newcommand{\pa}{\partial}
\newcommand{\al}{\alpha}
\newcommand{\ld}{\ldots}
\newcommand{\noa}{\noalign{\medskip}}

\begin{abstract}
Recall that a vector field on an $n$-dimensional differentiable manifold
$M$ is a mapping $X$ defined on $M$ with values in the tangent bundle $TM$
that assigns to each point $x\in M$ a vector $X(x)$ in the tangent space
$T_x M$. A vector field may be interpreted alternatively as the right-hand
side of an autonomous system of first-order ordinary differential equations,
i.e., a flow. Now we show that any flow can be enveloped by a conservative
dynamics using a semi-Riemann metric $g$ on $M$. This kind of dynamics was
called {\it geometric dynamics} [7]-[9].

The given vector field, the initial semi-Riemann metric, the Levi-Civita
connection, and an associated (1,1)-tensor field are used to build a new
geometric structure (e.g., semi-Riemann-Jacobi, semi-Riemann-Jacobi-Lagrange,
semi-Finsler-Jacobi, etc) on the manifold $M$ ensuring that all the trajectories
of a geometric dynamics are pregeodesics (Lorentz-Udri\c ste world-force law).
Implicitly, we solved a problem rised first by Poincar\'e: find a suitable
geometric structure that converts the trajectories of a given vector field
into geodesics.
\end{abstract}



\section{Geometric Dynamics}

Let $M$ be an $n$-dimensional differentiable manifold. A $C^\ty$ vector field
$X$ on $M$ defines the flow
$$
{dx \ov dt} = X(x). \leqno (1)
$$

A semi-Riemann metric $g$ on the manifold $M$ is a $C^\ty$ symmetric tensor
field of type (0,2) which assigns to each point $x\in M$ a nondegenerate
inner product $g(x)$ on the tangent space $T_x M$ of signature $(r,s)$. The
pair $(M,g)$ is called a {\it semi-Riemann manifold}.

The vector field $X$ and the semi-Riemann metric $g$ determine the {\it energy}
$$
f: M \to R, \quad f = {1 \ov 2} g(X,X).
$$
The vector field (flow) $X$ on $(M,g)$ is called:

1) {\it timelike}, if $f < 0$;

2) {\it nonspacelike or causal}, if $f \le 0$;

3) {\it null or lightlike}, if $f=0$;

4) {\it spacelike}, if $f>0$.

Let $X$ be a nonwhere vanishing vector field, of everywhere constant energy.
Upon rescaling, it may be supposed that $f \in \left\{ -1, 0, 1 \right\}$.
Generally, $\cal E$ is the set of zeros of the vector field
$X$. If $\cal E \ne \emptyset$, then this rescaling is possible only on
$M \setminus \cal E$.

Let $\na$ be the Levi-Civita connection of $(M,g)$. Using the operator
$\di{\na \ov dt}$  (covariant differentiation along a solution) we obtain the
prolongation
$$
{\na \ov dt} {dx \ov dt} = \na_{dx \ov dt} X \leqno (2)
$$
of the differential system (1) or of any perturbation of the system (1)
obtained adding to the second member $X$ a parallel vector field $Y$ with
respect to the covariant derivative $\na$. The prolongation by derivation
represents the general dynamics of the flow. The vector field $Y$ can be
used to illustrate a progression from stable to unstable flows, or converse.

The vector field $X$, the metric $g$, and the connection $\na$ determine the
external (1,1)-tensor field
$$
F = \na X - g^{-1} \otimes g (\na X),
$$
$$
F_j{}^i = \na_j X^i - g^{ih} g_{kj} \na_h X^k, \quad i,j,h,k = 1, \ld, n,
$$
which characterizes the {\it helicity} of vector field (flow) $X$.

First we write the differential system (2) in the equivalent form
$$
{\na \ov dt}{dx \ov dt} = g^{-1} \otimes g(\na X) \left( {dx \ov dt} \right)
+ F \left( {dx \ov dt} \right). \leqno (2')
$$

Successively we modify the differential system (2$'$) as follows:
$$
{\na \ov dt}  {dx \ov dt} = g^{-1} \otimes g(\na X)(X) +
F\left( {dx \ov dt} \right), \leqno (3)
$$
$$
{\na \ov dt} {dx \ov dt} = g^{-1} \otimes g (\na X) \left( {dx \ov dt} \right)
+ F(X), \leqno (4)
$$
$$
{\na \ov dt} {dx \ov dt} = g^{-1} \otimes g(\na X)(X) + F(X). \leqno (5)
$$
Obviously, the second order systems (3), (4), (5) are prolongations of the
first order system (1). Each of them is connected either to the dynamics
of the field $X$ or to the dynamics of a particle which is sensitive to the
vector field  $X$. Since
$$
g^{-1} \otimes g (\na X)(X) = grad\: f,
$$
we shall show that the prolongation (3) describes a conservative dynamics
of the vector field $X$ or of a particle which is
sensitive to the vector field $X$. The physical phenomenon produced by (4)
or (5) was not yet studied.

{\bf Theorem}. 1) {\it If $F=0$, then the kinematic system (1) prolonges to
a potential dynamical system with $n$ degrees of freedom, namely}
$$
{\na \ov dt} {dx \ov dt} = grad\: f. \leqno (3')
$$

2) {\it If $F \ne 0$, then the kinematic system (1) prolonges to a non-potential
dynamical system with $n$ degrees of freedom, namely}
$$
{\na \ov dt} {dx \ov dt} = grad\: f +
F\left( {dx \ov dt} \right). \leqno (3'')
$$

{\bf Corollary}. {\it If the metric $g$ is chosen such that
$f \in \left\{ -1, 0, 1 \right\}$ on $M \setminus \cal E$,
then the flow generated by $X$ is global and the dynamical systems (3$\:'$),
(3$\:''$) are reduced to}
$$
{\na^2 x \ov dt^2} = 0, \; {\na^2 x \ov dt^2} =
F\left( {dx \ov dt} \right).
$$

Let us show that the dynamical systems (3$\:'$) and (3$\:''$) are conservative.
To simplify the exposition we identity the tangent bundle $TM$ with the
cotangent bundle $T^*M$ using the semi-Riemann metric $g$.

{\bf Theorem}. 1) {\it The trajectories of the dynamical system (3$\:'$) are the
extremals of the Lagrangian}
$$
L = {1 \ov 2} g\left( {dx \ov dt}, {dx \ov dt} \right) + f(x).
$$

2) {\it The trajectories of the dynamical system (3$\:''$) are the extremals of
the Lagrangian}
$$
L = {1 \ov 2} g\left( {dx \ov dt} - X, \; {dx \ov dt} - X\right) =
{1 \ov 2} g\left( {dx \ov dt}, {dx \ov dt} \right)
- g \left(X, \di{dx \ov dt} \right) + f(x).
$$

3) {\it The dynamical systems (3$ \:'$) and (3$ \:''$) are conservative, the
Hamiltonian being the same for both cases, namely}
$$
H = {1 \ov 2} g \left( {dx \ov dt}, {dx \ov dt} \right) - f(x).
$$

The restriction of the Hamiltonian $H$ to the flow of the vector field $X$
is zero. Obviously the values of the Hamiltonian $H$ can be positive, negative or zero,
even if the metric $g$ is a Riemannian metric; therefore, just in this case,
there exist boundary-value problems associated to the differential system (3),
having three solutions (for example, the first corresponding to constant total energy
$H<0$, the second for $H=0$, and the third for $H>0$).

For the next theorem we recall that a {\it pregeodesic} is a smooth curve which
may be reparametrized to be a geodesic.

{\bf Theorem (Lorentz-Udri\c ste World-Force Law)}. 1) {\it Every non-cons-tant
trajectory of the dynamical system (3$ \: '$), which corresponds to a constant
value $H_0$ of the Hamiltonian, is a pregeodesic of the semi-Riemann-Jacobi
manifold}
$$
(M \setminus \cal E, \; \bar g = (H_0 + f) g).
$$

2) {\it Let $g_{ij}$ be the local components of the metric $g$ and let
$\Gamma^i_{jk}$, $i,j,k = 1, \ld, n$ be the local components of the connection
$\na$. Every non-constant trajectory of the dynamical system (3$\;''$), which
corresponds to a constant value $H_0$ of the Hamiltonian, is a horizontal
pregeodesic of the semi-Riemann-Jacobi-Lagrange manifold}
$$
(M \setminus \cal E, \; \bar g = (H_0 + f)g, \; N_j{}^i =
\Gamma^i_{jk} y^k - F_j{}^i, \quad i,j,k = 1, \ld, n).
$$

{\bf Corollary}. {\it If the metric g is chosen such that $f\in \{-1,0,1\}$
on $M \setminus \cal E$ and if we denote
$$
\al^2 = g \left( {dx \ov dt}, {dx \ov dt} \right), \quad
\beta = k^{1/2} g \left(X, {dx \ov dt} \right),
$$
then every trajectory of the dynamical system (3$\:''$), with $\al^2 k^{1/2}
+ \beta = 0$, is a pregeodesic of the semi-Finsler-Jacobi manifold}
$$
(M \setminus \cal E, \; L = F^2_w = w \al \beta, \quad w = \hbox{constant}).
$$

\section{Hamiltonian Structures on the Tangent Bundle}

Let $N$ be a $2n$-dimensional manifold. A nondegenerate and closed 2-form
$\wedge$ on $N$ is called {\it symplectic form}. A manifold $N$ with a given
symplectic form is called a {\it phase space}.

Let $(N, \wedge)$ be a phase space, and $H: N \to R$ be a $C^\ty$ real function.
We define the Hamilton gradient $X_H$ as being the vector field which satisfies
$$
\wedge_p (X_H (x), v) = dH (x)(v), \quad \forall v \in T_x N,
$$
and the Hamilton equations as
$$
{dx \ov dt} = X_H (x).
$$

Let $(M,g)$ be a semi-Riemann manifold with $n$ dimensions. Let $X$ be a
$C^\ty$ vector field on $M$, and $\omega = g \circ F$ the 2-form associated
to the tensor field $F = \na X - g^{-1} \otimes g(\na X)$ via the metric $g$.

The tangent bundle is usually endowed with the Sasaki metric $G$ created by
$g$. If $(x^i,y^i)$ are the coordinates of the point $(x,y) \in TM$ and
$\Gamma^i_{jk}$ are the components of the connection induced by $g_{ij}$,
then
$$
\left( {\pa \ov \pa x^i} - \Gamma^h_{ij} y^j {\pa \ov \pa y^h}, {\pa \ov \pa y^i}
\right), \quad (dx^j, \; \delta y^j = dy^j + \Gamma^j_{hk} y^k dx^h)
$$ are dual frames. Also the metric of Sasaki transcribes
$$
G = g_{ij} dx^i \otimes dx^j + g_{ij} \delta y^i \otimes \delta y^j.
$$

{\bf Theorem}. {\it The dynamical system (3$\:'$) lifts to $TM$ as a Hamilton
dynamical system with respect to the Hamiltonian
$$
H = {1 \ov 2} g \left( {dx \ov dt}, {dx \ov dt} \right) - f(x)
$$
and the symplectic 2-form}
$$
\Omega_1 = g_{ij} dx^i \wedge \delta y^j.
$$

{\bf Hint}. $\eta_1 = g_{ij} y^i dx^j$, and $d\eta_1 = -\Omega_1$.

{\bf Theorem}. {\it The dynamical system (3$\:''$) lifts to $TM$ as a Hamilton
dynamical system with respect to the Hamiltonian
$$
H = {1 \ov 2} g \left( {dx \ov dt}, {dx \ov dt} \right) - f(x)
$$
and the symplectic 2-form}
$$
\Omega_2 = {1 \ov 2} \omega_{ij} dx^i \wedge dx^j + g_{ij} dx^i \wedge
\delta y^j.
$$

{\bf Hint}. $\eta_2 = -g_{ij} X^i dx^j + g_{ij} y^i dx^j$, and
$d\eta_2 = -\Omega_2$.

\medskip

{\bf Pendulum Geometric Dynamics}. We use the Riemannian manifold
$(R^2, \delta_{ij})$. The small oscillations of a plane pendulum are described
as solutions of the differential system (plane pendulum flow)
$$
{dx_1 \ov dt} = -x_2, \quad {dx_2 \ov dt} = x_1. \leqno (6)
$$
In this case $x_1 (t) = 0$, $x_2 (t) = 0$, $t \in R$ is the equilibrium
point and $x_1 (t) = c_1 \cos t + c_2 \sin t$, $x_2 (t) = c_1 \sin t - c_2 \cos t$,
$t \in R$ is the general solution (family of circles with same centre).

Let
$$
X = (X_1, X_2), \; X_1 (x_1, x_2) = -x_2, \quad X_2 (x_1, x_2) = x_1,
$$
$$
f(x_1, x_2) = \di{1 \ov 2} (x^2_1 + x^2_2), \quad \hbox{rot}\: X = (0,0,2),
\; div \: X = 0.
$$

The pendulum flow conserves the areas.
The prolongation by derivation of the kinematic system (6) is
$$
{d^2x_i \ov dt^2} = \sum_j {\pa X_i \ov \pa x_j} {dx_j \ov dt}, \quad
i,j = 1,2
$$
or
$$
{d^2x_1 \ov dt^2} = -{dx_2 \ov dt}, \quad {d^2x_2 \ov dt^2} = {dx_1 \ov dt}.
\leqno (7)
$$
This prolongation admits the general solution
$$
\begin{array}{l}
x_1 (t) = a_1 \cos t + a_2 \sin t + h \\ \noa
x_2 (t) = a_1 \sin t - a_2 \cos t + k, \; t \in R. \end{array}
$$
(family of circles).

The pendulum geometric dynamics is described by
$$
{d^2x_i \ov dt^2} = {\pa f \ov \pa x_i} + \sum_j
\left( {\pa X_i \ov \pa x_j} - {\pa X_j \ov \pa x_i}\right) {dx_j \ov dt},
\quad i,j = 1,2
$$
or
$$
{d^2x_1 \ov dt^2} = x_1 - 2{dx_2 \ov dt}, \quad
{d^2x_2 \ov dt^2} = x_2 + 2{dx_1 \ov dt}, \leqno (8)
$$
with the general solution
$$
x_1 (t) = b_1 \cos t + b_2 \sin t + b_3 t\cos t + b_4 t \sin t
$$
$$
x_2 (t) = b_1 \sin t - b_2 \cos t + b_3 t\sin t - b_4 t \cos t, \; t \in R
$$
(family of spirals).

Using
$$
\begin{array}{l}
L = \di{1 \ov 2} \left[ \left( \di{dx_1 \ov dt}\right)^2 +
\left( \di{dx_2 \ov dt}\right)^2\right] + x_2 \di{dx_1 \ov dt} - x_1
\di{dx_2 \ov dt} + f \\ \noa
H = \di{1 \ov 2} \left[ \left( \di{dx_1 \ov dt}\right)^2 +
\left( \di{dx_2 \ov dt}\right)^2\right] - f \\ \noa
g_{ij} = (H + f) \delta_{ij}, \\ \noa
N_j{}^i = - F_j{}^i = - \delta^{ih} F_{jh}, \quad
F_{ij} = \di{\pa X_j \ov \pa x_i} - \di{\pa X_i \ov \pa x_j}, \; i,j,h =1,2,
\end{array}
$$
the solutions of the differential system (8) are horizontal pregeodesics of
the Riemann-Jacobi-Lagrange manifold
$$
(R^2 \setminus \{0\}, \; g_{ij}, N_j{}^i ).
$$

{\bf Lorenz Geometric Dynamics}. We use the Riemannian manifold $(R^3,
\delta_{ij})$. The Lorenz flow is a first dissipative model with chaotic
behaviour discovered in numerical experiment. Its state equations are
$$
\begin{array}{l}
\di{dx_1 \ov dt} = - \sigma x_1 + \sigma x_2 \\ \noa
\di{dx_2 \ov dt} = - x_1 x_3 + r x_1 - x_2 \\ \noa
\di{dx_3 \ov dt} =  x_1 x_2 - b x_3, \end{array} \leqno (9)
$$
where $\sigma, r, b$ are real parameters. Usually $\sigma, b$ are kept fixed
whereas $r$ is varied. At
$$
r > r_0 = {\si(\si + b+ 3) \ov \si - b-1}
$$
chaotic behaviour is observed. With $\si = 10, \; b=\di{8 \ov 3}$, the
preceding inequality yields $r_0 = 24, 7368$. If $\si \ne 0$ and
$b(r-1)>0$, then the equilibrium points of the Lorenz flow are
$$
x=0, \; y=0, \; z=0;
$$
$$
x = \pm \sqrt{b(r-1)}, \; y = \pm \sqrt{b(r-1)}, \; z= r-1.
$$

Let
$$
X = (X_1, X_2, X_3), \; X_1 (x_1, x_2, x_3) = -\si x_1 + \si x_2, \;
$$
$$
X_2 (x_1, x_2, x_3) = -x_1x_3 + rx_1 - x_2,\;
X_3(x_1,x_2,x_3) = x_1x_2 - bx_3,
$$
$$
f = {1 \ov 2} [(-\sigma x_1 + \si x_2)^2 + (-x_1x_3 + rx_1 - x_2)^2 +
(x_1x_2 - bx_3)^2],
$$
$$
rot X = (2x_1, \; -x_2, \; r-x_3 - \si).
$$

The Lorenz geometric dynamics is described by
$$
{d^2x_i \ov dt^2} = {\pa f \ov \pa x_i} + \sum_j
\left( {\pa X_i \ov \pa x_j} - {\pa X_j \ov \pa x_i} \right) {dx_j \ov dt},\;
i,j = 1,2,3
$$
or
$$
\begin{array}{l}
\di{d^2 x_1 \ov dt^2} = \di{\pa f \ov \pa x_1} + (\si + x_3 - r) \di{dx_2 \ov dt}
- x_2 \di{dx_3 \ov dt} \\ \noa
\di{d^2 x_2 \ov dt^2} = \di{\pa f \ov \pa x_2} + (r - x_3 - \si) \di{dx_1 \ov dt}
- 2x_1 \di{dx_3 \ov dt} \\ \noa
\di{d^2 x_3 \ov dt^2} = \di{\pa f \ov \pa x_3} + x_2 \di{dx_1 \ov dt} +
2x_1 \di{dx_2 \ov dt}. \end{array} \leqno (10)
$$

Using
$$
\begin{array}{l}
L = \di{1 \ov 2} \di\sum^3_{i=1} \left( {dx_i \ov dt}\right)^2 -
\di\sum^3_{i=1} X_i \di{dx_i \ov dt} + f \\ \noa
H = \di{1 \ov 2} \di\sum^3_{i=1} \left( \di{dx_i \ov dt}\right)^2 - f \\ \noa
g_{ij} = (H + f) \delta_{ij} \\ \noa
N_j{}^i = - F_j{}^i = - \delta^{ih}F_{jh}, \;
F_{ij} = \di{\pa X_j \ov \pa x_i} - \di{\pa X_i \ov \pa x_j}, \;
i,j,h = 1,2,3, \end{array}
$$
the solutions of the differential system (10) are horizontal pregeodesics
of the Riemann-Jacobi-Lagrange manifold
$$
(R^3 \setminus \cal E, \; g_{ij}, \; N_j{}^i),
$$
where $\cal E$ is the set of equilibrium points.

{\bf ABC Geometric Dynamics}. We use the Riemannian manifold $(R^3, \delta_{ij})$.
One examples of a fluid velocity that contains exponential stretching and
hence instability is the ABC flow,
$$
\left\{ \begin{array}{l}
\di{dx_1 \ov dt} = A\sin x_3 + C \cos x_2 \\ \noa
\di{dx_2 \ov dt} = B\sin x_1 + A \cos x_3 \\ \noa
\di{dx_3 \ov dt} = C\sin x_2 + B \cos x_1. \end{array} \right. \leqno (11)
$$

This flow is named after the three mathematicians Arnold, Beltrami and
Childress, who have contributed much to our understanding and appreciation of
classes of "chaotic" flows of which the present one is an example. For nonzero
values of the constants $A, B, C$ the preceding system is not globally integrable.
The topology of the flow lines is very complicated and can only be investigated
numerically to reveal regions of chaotic behaviour. The $ABC$ flow conserves
the volumes since the $ABC$ field is solenoidal.

The $ABC$ geometric dynamics is described by
$$
{d^2x^i \ov dt^2} = {\pa f \ov \pa x^i} + \sum_j \left( {\pa X_i \ov \pa x_j}
- {\pa X_j \ov \pa x_i} \right) {dx_j \ov dt}, \; i,j = 1,2,3.
$$
Since
$$
f = {1 \ov 2} (A+B+C+ 2AC\sin x_3 \cos x_2 + 2BA \sin x_1 \cos x_3 +
2 CB \sin x_2 \cos x_1)
$$
$$
rot X = X,
$$
the $ABC$ geometric dynamics is given by the differential system (12):
$$
\begin{array}{lcl}
\di{d^2 x_1 \ov dt^2} &=& AB \cos x_1 \cos x_3 - BC \sin x_1 \sin x_2 -
(B \cos x_1 + C \sin x_2) \di{dx_2 \ov dt} + \\ \noa
&+& (B \sin x_1 + A \cos x_3) \di{dx_3 \ov dt} \end{array}
$$
$$
\begin{array}{lcl}
\di{d^2 x_2 \ov dt^2} &=& - AC \sin x_2 \sin x_3 + BC \cos x_1 \cos x_2 +
(B \cos x_1 + C \sin x_2) \di{dx_1 \ov dt} - \\ \noa
&-& (A \sin x_3 + C \cos x_2) \di{dx_3 \ov dt} \end{array}
$$
$$
\begin{array}{lcl}
\di{d^2 x_3 \ov dt^2} &=& AC \cos x_3 \cos x_2 - BA \sin x_1 \sin x_3 -
(B \sin x_1 + A \cos x_3) \di{dx_1 \ov dt} + \\ \noa
&+& (C \cos x_2 + A \sin x_3) \di{dx_2 \ov dt}. \end{array}
$$

Using
$$
\begin{array}{l}
L = \di{1 \ov 2} \di\sum^3_{i=1} \left( \di{dx_i \ov dt} \right)^2 -
\di\sum^3_{i=1} X_i \di{dx_i \ov dt} + f \\ \noa
H = \di{1 \ov 2} \di\sum^3_{i=1} \left( \di{dx_i \ov dt} \right)^2 -f \\ \noa
g_{ij} = (H + f) \delta_{ij} \\ \noa
N_j{}^i = - F_j{}^i = - \delta^{ih}F_{jh}, \; F_{ij} = \di{\pa X_j \ov \pa x_i}
- \di{\pa X_i \ov \pa x_j}, \; i,j,h = 1,2,3 \end{array}
$$
the solutions of the differential system (12) are horizontal pregeodesics of
the Riemann-Jacobi-Lagrange manifold
$$
(R^3 \setminus \cal E, g_{ij}, N_j{}^i),
$$
where $\cal E$ is the set of equilibrium points which is included in the surface
of equation
$$
\sin x_1 \sin x_2 \sin x_3 + \cos x_1 \cos x_2 \cos x_3 = 0.
$$


\begin{thebibliography}{19}


\bibitem{[1]} R.G.Beil, {\it Comparison of unified field theories}, Tensor
N.S., 56, 2(1995), 175-183.

\bibitem{[2]} R.G.Beil, {\it Notes on a new Finsler metric function}, Balkan
Journal of Geometry and Its Applications, 2, 1(1997), 1-6.

\bibitem{[3]} F.M.Crampin, {\it A linear connection associated with any
second order differential equation field}, In: Eds. L.Tamassy, J.Szenthe,
New Developments in Geometry, 77-86, Kluwer Academic Publishers, 1996.

\bibitem{[4]} S.Friedlander, V.Yudovich, {\it Instabilities in
fluid motions}, Notices of the AMS, 46, 11(1999), 1358-1367.

\bibitem{[5]} V.Ob\u adeanu, C.Vernic, {\it Systemes dynamiques sur des
espace de Riemann}, Balkan Journal of Geometry and Its Applications, 2 1(1997),
73-82.

\bibitem{[6]} A.Udri\c ste, C.Udri\c ste, {\it Dynamics induced by a magnetic
field}, Conference on Differential Geometry, Budapest, July 27-30, 1996.
In: Ed.J.Szenthe, New Developments in Differential Geometry, Kluwer Academic
Publishers (1996), 429-442.

\bibitem{[7]} C.Udri\c ste, A.Udri\c ste, {\it Electromagnetic dynamical
systems}, Balkan Journal of Geometry and Its Applications, 2,1(1997),
129--140.

\bibitem{[8]} C.Udri\c ste, {\it Geometric dynamics}, Second Conference of
Balkan Society of Geometers, Aristotle University of Thessaloniki, June
23-26, 1998; Southeast Asian Bulletin of Mathematics 24, 1(2000), 1-11.

\bibitem{[9]} C.Udri\c ste, {\it Geometric dynamics}, Kluwer Academic
Publishers, 2000.

\bibitem{[10]} C.Udri\c ste, {\it Nonclassical Lagrangian Dynamics and Potential
Maps}, Proceedings of the Conference on Mathematics in Honour of Professor
Radu Ro\c sca at the Occasion of his Ninetieth Birthday, Katholieke University
Brussel, Katolieke University Leuven, Belgium, Dec. 11-16, 1999.

\bibitem{[11]} C.Udri\c ste, {\it Dynamics induced by second-order objects},
BSG Proceedings 4, Global Analysis, Differential Geometry, Lie Algebras,
Editor Grigorios Tsagas, pp.161-168, Geometry Balkan Press, 2000.

\bibitem{[12]} C.Udri\c ste, {\it Solutions of DEs and PDEs as potential
maps using first order Lagrangians}, Centennial Vr\u anceanu, Romanian
Academy, University of Bucharest, June 30 - July 4, 2000.

\bibitem{[13]} C.Udri\c ste, M.Neagu, {\it Geometrical interpretation of
solutions of certain PDEs}, Balkan Journal of Geometry and Its Applications,
4, 1(1999), 138-145.

\end{thebibliography}
\end{document}